\documentclass[a4paper,12pt,legno]{article}

\usepackage{amsmath}
\usepackage{amsfonts}
\usepackage{amsthm}
\usepackage{mathrsfs}
\usepackage{amssymb}

\newtheorem{prop}{PROPOSITION}
\newtheorem{cor}{COROLLARY}

\theoremstyle{remark}
\newtheorem{rem}{\textbf{Remark}}

\usepackage{graphicx}
\usepackage[T1]{fontenc}
\usepackage[polish]{babel}
\usepackage[cp1250]{inputenc}

\newcommand{\Rset}{\mathbb{R}}

\begin{document}

\title{Remarks on factorization property of some stochastic integrals\footnote{Research funded by Narodowe Centrum Nauki (NCN)
Dec2011/01/B/ST1/01257}}

\author{Zbigniew J. Jurek (University of Wroc\l aw\footnote{Part of this work
was done when Author was visiting
Indiana University, Bloomington,\qquad \qquad \qquad \qquad Indiana,
USA in 2013.}\,\,)}

\date{June 22, 2013}
\maketitle

\begin{quote}\textbf{Abstract.} In the paper Sato (2006) there are
introduced two families of improper random integrals and the
corresponding two convolution semigroups of infinitely divisible
laws on $\Rset^d$. Theorem 3.1 gives a relation (a factorization
property) between those two integrals. Here, using \emph{the random
integral mappings} $I^{h,r}_{(a,b]}$ (cf. the survey article Jurek
(2011)), we give a simpler proof that is  also valid for measures on
Banach spaces. Furthermore, using our technique we establish yet
other relations between those two families of improper stochastic
integrals.



\emph{Abbreviated title: Factorization of some random integrals}

\end{quote}

\newpage
In the last few decades there have appeared  many papers on random
integral representations of convolution subsemigroups of the
(master) semigroup, ID, of all infinitely probability distributions.
Jurek-Vervaat(1983) on the class, L, of selfdecomposable measures
seems to be one of the first in that area. For more references cf.
Jurek (2011), Sato (2006) and Maejima, Perez-Abreu and Sato (2012).
Some of the subsemigroups were introduced via the random integrals
while the others were described by transformations of the L\'evy
(spectral) measures of some infinitely divisible distributions. The
latter approach was presented already in Jurek (1990) and the
resulting measures were called there as\emph{ $\lambda$-mixtures}.
Most of that research was done in Euclidean spaces but we have also
techniques and proofs that are applicable in any infinite
dimensional separable Banach space.

\medskip
In this note using random integral technique we provide shorter and
simpler proofs of the factorization property of the two transforms
(integral operators) introduced in Sato (2006). It seems that the
general random integral method is more useful than  considerations
of some specific cases.

\medskip
\textbf{1.}  For an interval $(a,b]$ in the positive half-line, two
deterministic functions $h$ (\emph{space change}) and $r$
(\emph{inner clock time change}),  and a L\'evy process $Y_{\nu}(t),
t\ge0$ on a real separable Banach space $E$, where $\nu\in ID$ is
the law of random variable $Y_{\nu}(1)$, we consider the following
mapping ( or the operator):
\[
\nu\longmapsto
I^{h,r}_{(a,b]}(\nu):=\mathcal{L}\big(\int_{(a,b]}h(t)\,dY_{\nu}(r(t))\big)
\ \ (\star)
\]
and $\mathcal{L}$ denotes the probability distribution of the random
(stochastic) integral. Random integrals $(\ast)$ are defined by formal integration
by parts formula, i.e.,
\begin{multline*}
\int_{(a,b]} h(t)dY_{\nu}(r(t)):=
\\ h(b)Y_{\nu}(r(b))-h(a)Y_{\nu}(r(a))- \int_{(a,b]}Y_{\nu}(r(t)-)dh(t) \in E,
\end{multline*}
cf. Jurek and Vervaat (1983) or Jurek and Mason (1993) for a
discussion on the above random integrals.

Improper mappings $I^{h,r}_{(a,\infty)}$ are defined as limits as $b
\to \infty$; similarly, as limits, are defined the improper random
integrals $I^{h,r}_{(a,b)}$ cf. Jurek (2011) (invited Section
Lecture at $10^{th}$ Vilnius Conference on Probability in 2010) or
Jurek (2012).

Recall here that the integral $I^{h, r}_{(a,b]}$ commute with each
other, that is,
\[
I^{h_1, r_1}_{(a_1,b_1]}\big(I^{h_2, r_2}_{(a_2,b_2]}(\mu)\big)=
I^{h_2, r_2}_{(a_2,b_2]}\big(I^{h_1, r_1}_{(a_1,b_1]}(\mu)\big),
\]
provided $\mu$ is in appropriate domains. It follows from the
L\'evy-Khintchine formula for characteristic functions of infinitely
divisible distributions; cf. for details Jurek (2012).

\medskip
\textbf{2.}   For $-\infty<\beta<\alpha<\infty$, let us define the
following two families of time change clocks:
\begin{multline}
r_{\alpha}(t):= \int_t^{\infty} u^{-\alpha-1}e^{-u}du, \ \
\qquad \mbox{for} \ \ \ 0<t<\infty;  \ \ \ \ \ \mbox{and} \\
r_{\beta, \alpha}(t):= (\Gamma(\alpha - \beta))^{-1} \int_s^1
(1-u)^{\alpha-\beta -1}\,u^{-\alpha-1}du, \  \ \qquad \mbox{for} \ \
\ 0<t<1\,.
\end{multline}
Sato (2006) used the implicitly given inverse functions
$r_{\alpha}^{-1}$ and $r_{\alpha,\beta}^{-1}$ to define two improper
random integrals. In our notations these were random integral
mappings $ I^{t,\, r_{\alpha}(t)}_{(0,\infty)} \ \mbox{and} \ \
I^{s,\,\, r_{\alpha,\beta}(s)}_{(0,1)}$. One of the main result is
the following factorizations of the above two mappings:
\begin{prop}
 For $-\infty<\beta<\alpha<\infty$ and infinitely divisible $\nu$, on a real separable Banach space,
such that the following integrals are well defined we have that
\begin{multline}
I^{t,\,\, \int_t^{\infty} u^{-\beta-1}e^{-u}du}_{(0,\infty)}
(I^{s,\,\, (\Gamma(\alpha - \beta))^{-1} \int_s^1
(1-u)^{\alpha-\beta -1}\,u^{-\alpha-1}du}_{(0,1)}(\nu))\\
=I^{t,\,\, \int_t^{\infty} u^{-\alpha-1}e^{-u}du}_{(0,\infty)}(\nu)
\qquad \qquad \qquad
\end{multline}
\end{prop}
\begin{rem}
(i) Above we keep the explicite form the inner clock time for an
easy reference and comparison.

(ii) For general questions related to domains of the above random
integrals we refer to Sato (2006) and Jurek (2012). However, from
Jurek (2012),Corollary 10, we infer that in (2) for $\nu$ we can
take stable measures with the exponent $p>\alpha$.

(iii) Also, the proof below is valid for any real separable infinite
dimensional Banach space -- not only for Euclidean space $\Rset^d$
as it is in Sato (2006).
\end{rem}

\medskip
\emph{Proof of Proposition 1.} As in Theorem 2, Section 4.2 in Jurek
(2012), let us define Borel  measures $\rho_i$ using the inner clock
time change from (1). Namely, let
\begin{equation}
\rho_1((c,d]):= \int_{(c,d]} u^{-\beta-1}e^{-u}du , \ \ \ (c,d]
\subset(0, \infty)
\end{equation}
and
\begin{equation}
\rho_2((c,d]):= (\Gamma (\alpha-\beta))^{-1}\int_{(c,d]}(1-u)^{\alpha-\beta -1} \,u^{-\alpha-1}du,
\ \ \ (c,d] \subset(0,1)
\end{equation}
Furthermore, let define the space change functions as follows
\begin{equation}
h_1(t):=t, \, \, t\in (0,\infty) \qquad \mbox{and} \qquad h_2(s):=s,
\, \, s\in (0,1)
\end{equation}
Finally, let
\begin{equation}
\boldsymbol{\rho}: = \rho_1 \times \rho_2  \ \mbox{and} \ \ \
\textbf{h}(t,s):=h_1 \otimes h_2 (t,s)=h_1(t)h_2(s) \ \mbox{(tensor
product)}
\end{equation}
Now observe that for the image measure $\textbf{h}\boldsymbol{\rho}$ and $u>0$ we have
\begin{multline*}
(\textbf{h}\boldsymbol{\rho})(x:x>u)=\int_0^{\infty}1_{(x:x>u)}(v)\textbf{h}\boldsymbol{\rho}(dv)
=\int_0^{\infty}\int_0^1\,1_{(x:x>u)}(s\cdot t)\rho_1(ds)\rho_2(dt)\\
=(\Gamma (\alpha-\beta))^{-1}\int_0^{\infty}\big( \int_0^1\,1_{(x:x>u)}(s\cdot t)(1-s)^{\alpha-\beta -1} \,s^{-\alpha-1}ds\big)t^{-\beta-1}e^{-t}dt  \ \ \mbox{(w:=st)}  \\
=(\Gamma (\alpha-\beta))^{-1} \int_0^{\infty}\big( \int_0^t\,1_{(x:x>u)}(w)(1-\frac{w}{t})^{\alpha-\beta -1} \,(\frac{w}{t})^{-\alpha-1}\frac{dw}{t}\big)t^{-\beta-1}e^{-t}dt  \\
=(\Gamma (\alpha-\beta))^{-1} \int_0^{\infty}\big( \int_0^t\,1_{(x:x>u)}(w)(t-w)^{\alpha-\beta -1} \,w^{-\alpha-1}dw\big)e^{-t}\,dt \ \ \mbox{(changing order)}\\
=(\Gamma (\alpha-\beta))^{-1}\int_0^{\infty}1_{(x:x>u)}(w) w^{-\alpha-1}\big( \int_w^{\infty}(t-w)^{\alpha-\beta -1}e^{-t}dt\big)dw  \\
=(\Gamma (\alpha-\beta))^{-1}\int_0^{\infty}1_{(x:x>u)}(w) w^{-\alpha-1} e^{-w}\big(\int_w^{\infty}(t-w)^{\alpha-\beta -1}e^{-(t-w)}dt\big)dw      \\
=\int_u^{\infty} w^{-\alpha-1}e^{-w}dw.
\end{multline*}
Hence and from Theorem 2 in Jurek (2012) we get the equality (2)
which completes the proof.

\begin{cor} (a) For $-\infty<\beta<\alpha<\infty$ and the inner clock changes
$r_{\alpha}$ and $r_{\beta,\,\alpha}$ given in (1) we have a
factorization
\[
I^{t,\, r_{\beta}(t)}_{(0,\infty)} \circ
I^{s,\,r_{\beta,\alpha}(s)}_{(0,1)}
=I^{t,\,r_{\alpha}(t)}_{(0,\infty)}\;
\]
(b) For $-\infty<\alpha_k<\alpha_{k-1}<\alpha_{k-2}<. . . <
\alpha_2<\alpha_1<\infty$ we have
\[
I^{t,\, r_{\alpha_k}(t)}_{(0,\infty)}\circ
I^{s,\,r_{\alpha_k,\alpha_{k-1}}(s)}_{(0,1)} \circ
I^{s,\,r_{\alpha_{k-1},\alpha_{k-2}}(s)}_{(0,1)}\circ...\circ
I^{s,\,r_{\alpha_2,\alpha_1}(s)}_{(0,1)} = I^{t,\,
r_{\alpha_1}(t)}_{(0,\infty)}\,,
\]
where $\circ$ denotes the composition of the random integral
mappings.
\end{cor}
Proofs follows from Proposition 1 by mathematical induction
argument.

\medskip
\textbf{3.}   The following factorization was predicted but not
proved in Sato (2006) in Comment 2 on p. 86. Here it is phrased in
the terms of our integral mappings $I^{h,r}_{(a,b]}$.
\begin{prop}
For $ -\infty <\gamma <\beta<\alpha<\infty$ and an infinitely
divisible $\nu$, on a real separable Banach space, such that the
following integral are well defined, we have that
\begin{equation}
I^{t,\,r_{\beta,\alpha}(t)}_{(0,1)}\big(I^{s,\, r_{\gamma,
\beta}(s)}_{(0,1)}(\nu)\big)= I^{s,\, r_{\gamma,
\beta}(s)}_{(0,1)}\big(I^{t,\,r_{\beta,\alpha}(t)}_{(0,1)}(\nu)\big)
=I^{u,\,r_{\gamma,\alpha}(u)}_{(0,1)}(\nu)
\end{equation}
\end{prop}
\emph{Proof of Proposition 2.} For later use let recall the relation
between the special functions beta and gamma. Namely, for $a>0,\
b>0$
\[
B(a,b):=\int_0^1(1-u)^{a-1}u^{b-1}du, \ \ \ \mbox{and} \ \  \
B(a,b)=\frac{\Gamma(a)\,\Gamma (b)}{\Gamma (a+b)}
\]

As in proof of Proposition 1 we use Theorem 2 from Jurek (2012).

For the L\'evy exponent of the ID measure on the right hand side in
(7) we have
\begin{multline*}
\Gamma(\alpha-\beta)\,\Gamma(\beta-\gamma)\,\int_0^1\big(\int_0^1\,\Phi(s\,t\,y)|dr_{\alpha,\beta}(t)|\big)|dr_{\beta,
\gamma}(t)|\\=
\int_0^1\big(\int_0^1\,\Phi(s\,t\,y)(1-t)^{\alpha-\beta-1}t^{-\alpha-1}dt\big)(1-s)^{\beta-\gamma-1}s^{-\beta-1}ds
\ \ \
(\mbox{put}\ st=:w)  \\
=
\int_0^1\big(\int_0^s\,\Phi(w\,y)(1-\frac{w}{s})^{\alpha-\beta-1}(\frac{w}{s})^{-\alpha-1}\frac{dw}{s}\big)
(1-s)^{\beta-\gamma-1}s^{-\beta-1}ds \ \ (\mbox{change order})\\
=\int_0^1\Phi(wy)w^{-\alpha-1}\big(\int_w^1(s-w)^{\alpha-\beta-1}(1-s)^{\beta-\gamma-1}ds\big)dw
\ \ \ (\mbox{put} \ \ 1-s=:z)\\
=\int_0^1\Phi(wy)w^{-\alpha-1}\big(\int_0^{1-w}(1-w-z)^{\alpha-\beta-1}z^{\beta-\gamma-1}dz\big)dw
\ \ (\mbox{put} \ (1-w)^{-1}z=:x) \\
=\int_0^1\Phi(wy)w^{-\alpha-1}(1-w)^{\alpha-\gamma-1}dw \,\,
\int_0^1(1-x)^{\alpha-\beta -1}x^{\beta-\gamma -1}dx \\
=
B(\alpha-\beta,\beta-\gamma)\,\Gamma(\alpha-\gamma)\,\int_0^1\Phi(wy)|dr_{\alpha,\gamma}(w)|
\\=
\Gamma(\alpha-\beta)\, \,
\Gamma(\beta-\gamma)\int_0^1\Phi(wy)|dr_{\alpha,\gamma}(w)|,
\end{multline*}
which proves identity (7) and Proposition 2.
\begin{cor}
For positive integer $k\ge2$ and reals $\alpha_i, \, i=1,2,...,k$
such that
$-\infty<\alpha_k<\alpha_{k-1}<...<\alpha_2<\alpha_1<\infty$ we have
\[
I^{t,\,\,r_{\alpha_2,\alpha_1}(t)}_{(0,1)}
\circ\,I^{t,\,\,r_{\alpha_3,\alpha_2}(t)}_{(0,1)}
\circ\,I^{t,\,\,r_{\alpha_4,\alpha_3}(t)}_{(0,1)} \circ ... \circ
I^{t,\,\,r_{\alpha_k,\alpha_{k-1}}(t)}_{(0,1)}=
I^{t,\,\,r_{\alpha_k,\alpha_1}(t)}_{(0,1)}
\]
where $\circ$ denotes the composition of the random integral
mappings.
\end{cor}
Its proof follows from Proposition 2 via the induction argument.

\medskip
Last but not least, from the few instances showed in this note, one
may expect that the images of measures through tensor product will
find more applications and may provide simpler proofs as well; cf.
Jurek (2012).

{}

\medskip
\medskip
\medskip
\noindent Institute of Mathematics, University of Wroc\l aw, Pl.
Grunwaldzki 2/4, 50-384 Wroc\l aw, Poland. [E-mail:
zjjurek@math.uni.wroc.pl]

\end{document}